# INVERSE PROBLEM OF THE SPECTRAL ANALYSIS FOR THE STURM-LIOUVILLE OPERATOR WITH NON-SEPARATED BOUNDARY CONDITIONS AND SPECTRAL PARAMETER IN THE BOUNDARY CONDITION

**I.M. Nabiev**

*Department of Applied Mathematics, Baku State University, Z. Khalilov 23, AZ1148, Baku, Azerbaijan;*
*Institute of Mathematics and Mechanics, Azerbaijan National Academy of Sciences, B. Vahabzadeh 9, AZ1141, Baku, Azerbaijan;*
*Department of Mathematics, Khazar University, Mahsati 11, AZ1096, Baku, Azerbaijan.*
*E-mail:* nabievim@yahoo.com

**Abstract.** This work deals with an inverse problem for the Sturm-Liouville operator with non-separated boundary conditions, one of which linearly depends on a spectral parameter. Uniqueness theorem is proved, solution algorithm is constructed and sufficient conditions for solvability of inverse problem are obtained. As spectral data, we only use the spectrum of one boundary value problem and some sequence of signs.

## 1. Introduction

Many issues of mathematical physics are reduced to inverse problems of spectral analysis for various differential operators, which require the recovery of operators from some of their given spectral data. These problems are often considered in mathematics and various branches of natural science and technical science. They occur, for example, in quantum mechanics when defining intra-atomic forces from the known energy levels, in radiotechnics when synthesizing the parameters of nonhomogeneous transmission lines, in elasticity theory when defining the size of cross sections of the beam from the given frequencies of its self-oscillations, etc. Note that the research in the field of inverse problems has been given impetus also by the fact that they allowed creating new methods for the calculation of radio waves, thermal radiation and vibration of metal structures.

The object most closely related to the inverse problems is the Sturm-Liouville operator. Classical results by G. Borg, N. Levinson, M.G. Krein, M.G. Gasymov, B.M. Levitan, V.A. Marchenko dedicated to the problems of recovery of regular Sturm-Liouville operator with separated boundary conditions are well-known (see [1-2] and references therein). Besides, inverse problems for this operator with discontinuity conditions have been considered in [3-6], etc. Some variations of inverse problems with non-separated boundary conditions have been fully solved in well-known works by I.V. Stankevich, V.A. Sadovnichii, V.A. Marchenko, V.A. Yurko, O.A. Plaksina, M.G. Gasymov, I.M. Guseinov, I.M. Nabiev, Y.T. Sultanaev, A.M. Akhtyamov, Y.L. Korotyaev, A.S. Makin (see [7-18]).

Inverse problems for differential operators with a spectral parameter in the equation and boundary conditions are playing an important role in many physical and technical applications [16, 19]. In such problems, the unknown coefficients of differential equation and boundary conditions are recovered either from the spectra of two or three boundary value problems with different separated or non-separated boundary conditions, or from two spectra and two or three eigenvalues, or from the spectra of two boundary value problems, some sequence of signs and some number [20-35].

In this work, we consider the inverse spectral problem of recovery of Sturm-Liouville operator with non-separated boundary conditions, one of which contains a spectral parameter. Uniqueness theorem is proved, solution algorithm is constructed and sufficient conditions for solvability of inverse problem are obtained. In contrast to the problems considered earlier, here only one spectrum and a sequence of signs are used as spectral data.

**2. Physical interpretation of the problem. Inverse problem statement**

Consider the boundary value problem generated on the interval $[0, \pi]$ by the Sturm-Liouville equation

$$-y'' + q(x)y = \lambda^2 y \qquad (1)$$

and the non-separated boundary conditions of the form

$$y'(0)+(\alpha\lambda+\beta)y(0)+\omega y(\pi)=0,$$
$$y'(\pi)+\gamma y(\pi)-\omega y(0)=0, \quad (2)$$

where $q(x)$ is a real function belonging to the space $L_2[0,\pi]$, $\lambda$ is a spectral parameter, and $\alpha, \beta, \gamma, \omega$ are real numbers with $\alpha \neq 0$, $\omega \neq 0$. We denote this problem by $P$.

Let's give a physical (mechanical) interpretation of this problem. In [19, p. 31], a real mechanical system with non-separated boundary conditions is described and a scheme for rod is presented. But it can be assumed that it is not a rod but a string, one end of which is connected to the other by a rope with spring, and, besides, elastically fixed at both ends. Differential equation (1) describes the vibrations of a string, which vibrates in a variable medium with elasticity coefficient $q(x)$. Function $y(x)$ represents transverse displacement of a string. If the spring on the rope possess rigidity $c_2$, the ends of the string are fixed by the springs of rigidities $c_0$ and $c_1$, respectively, and there is a friction on the left end (friction is proportional to velocity, which is a time derivative), then the boundary conditions have the following form:

$$y'(0)=(\alpha\lambda+c_0+c_2)y(0)+c_2 y(\pi),$$
$$-y'(\pi)=(c_1+c_2)y(\pi)+c_2 y(0),$$

where $\alpha$ is a friction coefficient on the left end. Note that, based on [16, p. 214-217], one can provide another physical interpretation for the problem $P$, related to electrical vibrations in the wire, closed by concentrated resistance, inductance and a tank.

Denote by $c(x,\lambda)$, $s(x,\lambda)$ a fundamental system of solutions to the equation (1), defined by the initial conditions $c(0,\lambda)=s'(0,\lambda)=1$, $c'(0,\lambda)=s(0,\lambda)=0$. Let

$$\eta(x,\lambda)=c'(x,\lambda)+\gamma c(x,\lambda), \quad \sigma(x,\lambda)=s'(x,\lambda)+\gamma s(x,\lambda).$$

It is easy to see that the characteristic function of the boundary value problem $P$ is

$$\delta(\lambda)=2\omega-\eta(\pi,\lambda)+\omega^2 s(\pi,\lambda)+(\alpha\lambda+\beta)\sigma(\pi,\lambda). \quad (3)$$

The zeros $\mu_k$, $k = \pm 0, \pm 1, \pm 2, \ldots$ of the function $\delta(\lambda)$ are the eigenvalues of the problem $P$. Denote $\sigma_n = \text{sign}\left[1 - |\omega s(\pi, \theta_n)|\right]$ ($n = \pm 1, \pm 2, \ldots$), where $\theta_n$'s are the zeros of the function $\sigma(\pi, \lambda)$, whose squares are the eigenvalues of the boundary value problem generated by the equation (1) and boundary conditions $y(0) = y'(\pi) + \gamma y(\pi) = 0$.

This work is dedicated to the inverse problem stated as follows.

**Inverse problem.** Given the spectral data $\{\mu_k\}$, $\{\sigma_n\}$ of boundary value problems, construct the coefficient $q(x)$ of the equation (1) and the parameters $\alpha, \beta, \gamma, \omega$ of the boundary conditions (2).

### 3. Uniqueness theorem. Solution algorithm

**Theorem 1.** *The boundary value problem $P$ can be uniquely recovered if its spectrum $\{\mu_k\}$ and the sequence of signs $\{\sigma_n\}$ are known.*

**Proof.** As is known [17], the following representations are true for the functions $c(\pi, \lambda)$, $c'(\pi, \lambda)$, $s(\pi, \lambda)$ and $s'(\pi, \lambda)$:

$$c(\pi, \lambda) = \cos \lambda \pi + Q \frac{\sin \lambda \pi}{\lambda} + \frac{f_1(\lambda)}{\lambda}, \quad c'(\pi, \lambda) = -\lambda \sin \lambda \pi + Q \cos \lambda \pi + f_2(\lambda),$$

$$s(\pi, \lambda) = \frac{\sin \lambda \pi}{\lambda} - Q \frac{\cos \lambda \pi}{\lambda^2} + \frac{f_3(\lambda)}{\lambda^2}, \quad s'(\pi, \lambda) = \cos \lambda \pi + Q \frac{\sin \lambda \pi}{\lambda} + \frac{f_4(\lambda)}{\lambda},$$

where $Q = \frac{1}{2}\int_0^\pi q(x)dx$, $f_m(\lambda)$ ($m = 1, 2, 3, 4$) is an entire function of exponential type not greater than $\pi$, square summable on the real axis. Using these representations and Paley-Wiener theorem, from (3) we obtain that the characteristic function of the boundary value problem $P$ has the form

$$\delta(\lambda) = 2\omega + \lambda(\sin \lambda \pi + \alpha \cos \lambda \pi) + (\gamma + Q)\alpha \sin \lambda \pi + (\beta - \gamma - Q)\cos \lambda \pi + f(\lambda), \qquad (4)$$

where $f(\lambda) = \int_{-\pi}^{\pi} \tilde{f}(t) e^{i\lambda t} dt$, $\tilde{f}(t) \in L_2[-\pi, \pi]$. Using (4), it is easy to prove that the zeros of the function $\delta(\lambda)$ obey the asymptotics

$$\mu_k = k + a + \frac{(-1)^{k+1} A - B}{k\pi} + \frac{\tau_k}{k} \qquad (5)$$

as $|k| \to \infty$, where $\{\tau_k\} \in l_2$,

$$a = -\frac{1}{\pi}\operatorname{arctg}\alpha, \quad A = \frac{2\omega}{\sqrt{1+\alpha^2}}, \tag{6}$$

$$B = \frac{\beta}{1+\alpha^2} - \gamma - Q. \tag{7}$$

Using (5) and (6), define the parameters $\alpha$ and $\omega$ by the formulas

$$\alpha = -\operatorname{tg}\pi a, \quad \omega = \frac{\pi\sqrt{1+\alpha^2}}{2}\lim_{k\to\infty} k(\mu_{2k+1} - \mu_{2k} - 1), \tag{8}$$

where $a = \lim_{k\to\infty}(\mu_k - k)$. Knowing the spectrum $\{\mu_k\}$ and the parameter $\alpha$, we can recover the function $\delta(\lambda)$ (an entire function of exponential type) in the form of infinite product as follows:

$$\delta(\lambda) = \pi\sqrt{1+\alpha^2}\,(\mu_{-0} - \lambda)(\mu_{+0} - \lambda)\prod_{\substack{k=-\infty \\ k\neq 0}}^{\infty} \frac{\mu_k - \lambda}{k}. \tag{9}$$

Construct the function

$$\sigma(\pi,\lambda) = \frac{\delta(\lambda) - \delta(-\lambda)}{2\alpha\lambda}. \tag{10}$$

From (10) we find the zeros $\theta_n$, $n = \pm 1, \pm 2, \ldots$ of the function $\sigma(\pi,\lambda)$, for which the following asymptotic formula holds as $|n| \to \infty$:

$$\theta_n = n - \frac{1}{2}\operatorname{sign} n + \frac{Q+\gamma}{n\pi} + \frac{\xi_n}{n}, \tag{11}$$

where $\{\xi_n\} \in l_2$. As (10) is an even function, we have $\theta_{-n} = -\theta_n$. Knowing $\mu_k$ and $\theta_n$, we can also recover the parameter $\beta$, because, by (5), (7) and (11),

$$\beta = 2\omega - 2\pi(1+\alpha^2)\lim_{k\to+\infty} k\left(\mu_{2k+1} - \theta_{2k+1} - a - \frac{1}{2}\right). \tag{12}$$

Consider the functions

$$v_+(\lambda) = -\eta(\pi,\lambda) + \omega^2 s(\pi,\lambda), \tag{13}$$

$$v_-(\lambda) = -\eta(\pi,\lambda) - \omega^2 s(\pi,\lambda). \tag{14}$$

In view of (3), the function $v_+(\lambda)$ is recovered by the formula

$$v_+(\lambda) = \delta(\lambda) - 2\omega - (\alpha\lambda + \beta)\sigma(\pi,\lambda). \tag{15}$$

Now let's show that, in addition to the spectrum $\{\mu_k\}$ (which, as shown above, makes it possible to uniquely recover $v_+(\lambda)$, $\sigma(\pi,\lambda)$, $\alpha$, $\beta$, $\omega$), it suffices to prescribe a sequence $\{\sigma_n\}$ in order to recover the function $v_-(\lambda)$, and, consequently, the function

$$s(\pi,\lambda) = \frac{1}{2\omega^2}[v_+(\lambda) - v_-(\lambda)]. \tag{16}$$

In fact, as $\theta_n$'s are the roots of the function $\sigma(\pi,\lambda)$, it follows from the relation

$$c(\pi,\lambda)\sigma(\pi,\lambda) - s(\pi,\lambda)\eta(\pi,\lambda) = 1$$

that

$$s(\pi,\theta_n)\eta(\pi,\theta_n) = -1. \tag{17}$$

By this equality, it is easy to obtain from (13), (14) that $v_-^2(\theta_n) - v_+^2(\theta_n) = -4\omega^2$. Therefore,

$$v_-(\theta_n) = \text{sign } v_-(\theta_n)\sqrt{v_+^2(\theta_n) - 4\omega^2}. \tag{18}$$

Taking into account the intermittency of the zeros of the functions $s(\pi,\lambda)$ and $\sigma(\pi,\lambda)$ and the representation of the function $s(\pi,\lambda)$, we have $\text{sign } s(\pi,\theta_n) = (-1)^{n+1}$. Then, by (14) and (17),

$$\text{sign } v_-(\theta_n) = \text{sign}\left[-\eta(\pi,\theta_n) - \omega^2 s(\pi,\theta_n)\right] =$$

$$= \text{sign}\left[\frac{1}{s(\pi,\theta_n)} - \omega^2 s(\pi,\theta_n)\right] = \text{sign } \frac{1 - [\omega s(\pi,\theta_n)]^2}{s(\pi,\theta_n)} = (-1)^{n+1}\sigma_n.$$

Substituting this into (18), we obtain

$$v_-(\theta_n) = (-1)^{n+1}\sigma_n\sqrt{v_+^2(\theta_n) - 4\omega^2}. \tag{19}$$

Let

$$g(\lambda) = v_+(\lambda) - v_-(\lambda) - 2\omega^2\frac{\sin\lambda\pi}{\lambda}. \tag{20}$$

Using the asymptotics (11), the relations (13) and (14), the representation of the function $s(\pi,\lambda)$, Lemma 1.4.3 of [17] and Theorem 28 of [36], it is easy to obtain that the function $g(\lambda)$ is uniquely determined by the sequences $\{\theta_n\}$, $\{\sigma_n\}$, $\{v_+(\theta_n)\}$ by the formula

$$g(\lambda) = 2\sigma(\pi,\lambda)\sum_{n=1}^{\infty}\frac{\theta_n g(\theta_n)}{(\lambda^2 - \theta_n^2)\dot\sigma(\pi,\theta_n)}, \tag{21}$$

where $g(\theta_n) = \dot{v}_+(\theta_n) - (-1)^{k+1}\sigma_n \sqrt{v_+^2(\theta_n) - 4\omega^2} - 2\omega^2 \frac{\sin\theta_n \pi}{\theta_n}$, and the point over the function means differentiation with respect to $\lambda$.

The uniqueness of the constructed function $g(\lambda)$ follows from the fact that the interpolation formula (21) performs a bijection between $l_2$ and the space of entire functions of exponential type not greater than $\pi$, square summable on the real axis. Therefore, the characteristic function $s(\pi, \lambda)$ of the boundary value problem generated by the equation (1) and Dirichlet boundary conditions

$$y(0) = y(\pi) = 0 \tag{22}$$

is recovered by the formula (16), where the function $v_-(\lambda)$ is defined from (20). The zeros $\lambda_n$, $n = \pm 1, \pm 2, \ldots$ of the function (16) (whose squares are the eigenvalues of the boundary value problem (1), (22)) satisfy the asymptotic formula

$$\lambda_n = n + \frac{Q}{n} + \frac{\eta_n}{n}, \quad \{\eta_n\} \in l_2. \tag{23}$$

From (11) and (23) it follows

$$\gamma = \pi \lim_{n \to \infty} k\left(\theta_n - \lambda_n + \frac{1}{2}\right). \tag{24}$$

Finally, the characteristic function of the problem generated by the equation (1) and the boundary conditions

$$y(0) = y'(\pi) = 0 \tag{25}$$

is recovered by means of $\sigma(\pi, \lambda)$, $s(\pi, \lambda)$ and $\gamma$ by the following formula:

$$s'(\pi, \lambda) = \sigma(\pi, \lambda) - \gamma s(\pi, \lambda). \tag{26}$$

As is known [1, 17], the coefficient $q(x)$ of the equation (1) is defined uniquely by the zeros $v_n$, $n = 1, 2, \ldots$ of this function and the sequence $\{\lambda_n\}$.

Thus, the boundary value problem $P$ is completely recovered by means of the given sequences $\{\mu_k\}$ and $\{\sigma_n\}$. Theorem is proved.

Based on the proof of Theorem 1, let's state the solution algorithm for an inverse problem.

**Algorithm.** Let the sequences $\{\mu_k\}$, $\{\sigma_n\}$, i.e. the spectral data of the boundary value problem $P$, be given.

1) Define the parameters $\alpha$ and $\omega$ of the boundary conditions (2) from (8).

2) Construct the function $\delta(\lambda)$ in the form of infinite product (9).

3) Recover the function $\sigma(\pi, \lambda)$ by (10) and find the zeros $\theta_n$ of this function.

4) Calculate the parameter $\beta$ by the formula (12).

5) Recover the function (13) by (15).

6) Find the values of the function (14) at the points $\theta_n$ by using (19).

7) Using $\sigma(\pi, \lambda)$ and $v_+(\theta_n)$, recover the function (20) by the interpolation formula (21).

8) Knowing $g(\lambda)$, define $v_-(\lambda)$ from (20).

9) Define the characteristic function $s(\pi, \lambda)$ of the boundary value problem (1), (22) by the formula (16).

10) Using the asymptotic formulas (11) and (23), find the parameter $\gamma$ by the formula (24).

11) Using $\sigma(\pi, \lambda)$, $s(\pi, \lambda)$ and $\gamma$, recover the characteristic function $s'(\pi, \lambda)$ of the problem (1), (25) by the formula (26).

12) Using the sequences $\{\lambda_n\}$ and $\{v_n\}$ of zeros of the functions $s(\pi, \lambda)$ and $s'(\pi, \lambda)$, respectively, construct the coefficient $q(x)$ of the equation (1) following the procedure described in [1, 17].

**4. Sufficient conditions for solvability of inverse problem**

**Theorem 2.** *In order for the sequences of real numbers $\{\mu_k\}$ $(k = \pm 0, \pm 1, \pm 2,...)$ and $\{\sigma_n\}$ $(\sigma_n = -1, 0, 1;\ n = \pm 1, \pm 2,...)$ to be spectral data of the boundary value problem of the form* $P$, *it is sufficient that the following conditions hold:*

1) *the asymptotic formula (5) is true, where* $A = 2\omega \cos \pi a$, $\omega$, $a$, $B$ *are real numbers*, $0 < |a| < \dfrac{1}{2}$, $\omega \neq 0$, $\{\tau_k\} \in l_2$;

2) $... \leq \theta_{-3} \leq \mu_{-2} \leq \theta_{-2} \leq \mu_{-1} \leq \theta_{-1} \leq \mu_{-0} < 0 < \mu_{+0} \leq \theta_1 \leq \mu_1 \leq \theta_2 \leq \mu_2 \leq \theta_3 \leq ...$,

*where $\theta_n$'s are the zeros of the function $\delta(\lambda) - \delta(-\lambda)$,*

$$\delta(\lambda) = \frac{\pi(\mu_{-0} - \lambda)(\mu_{+0} - \lambda)}{\cos \pi a} \prod_{\substack{k=-\infty \\ k \neq 0}}^{\infty} \frac{\mu_k - \lambda}{k}, \qquad (27)$$

*with $\theta_n \neq \theta_m$ for $n \neq m$;*

*3) the inequality $b_n \stackrel{def}{=} |\delta(\theta_n) - 2\omega| - 2|\omega| \geq 0$ holds;*

*4) $\sigma_n$ is equal to zero if $b_n = 0$, and to 1 or -1 if $b_n > 0$; besides, there exists $N > 0$ such that $\sigma_n = 1$ for all $|n| \geq N$.*

**Proof.** Denote $\alpha = -\text{tg}\,\pi a$. As $\dfrac{1}{\cos \pi a} = \sqrt{1+\alpha^2}$ and the numbers $\mu_k$ obey the asymptotics (5), by Lemma 1.3 of [37], for the function (27) we have

$$\delta(\lambda) = \sqrt{1+\alpha^2}\,[A + \lambda \sin \pi(\lambda - a) + B \cos \pi(\lambda - a) + f(\lambda - a)],$$

where $f(\lambda) = M \sin \lambda \pi + \int_{-\pi}^{\pi} g(t)e^{i\lambda t}\,dt$, $g(t) \in L_2[-\pi, \pi]$, and $M$ is some constant. Hence, taking into account $A = 2\omega \cos \pi a$ and $\alpha = -\text{tg}\,\pi a$, we have

$$\delta(\lambda) = 2\omega + \lambda(\sin \lambda \pi + \alpha \cos \lambda \pi) + B_1 \cos \lambda \pi +$$

$$+ B_2 \sin \lambda \pi + \int_{-\pi}^{\pi} g_1(t)e^{i\lambda t}\,dt, \qquad (28)$$

where $B_1 = B + \alpha M$, $B_2 = M - \alpha B$, $g_1(t) = \sqrt{1+\alpha^2}\,g(t)e^{-i a t}$. Consequently, we get the following representation for the function $\sigma(\lambda) = \dfrac{\delta(\lambda) - \delta(-\lambda)}{2\alpha\lambda}$:

$$\sigma(\lambda) = \cos \lambda \pi - B_3 \pi \frac{\sin \lambda \pi}{\lambda} + \int_0^{\pi} \tilde{g}(t)\frac{\sin \lambda t}{\lambda}\,dt, \qquad (29)$$

where $B_3 = -\dfrac{B_2}{\pi \alpha}$, $\tilde{g}(t) \in L_2[0, \pi]$. By Lemma 3.4.2 of [17] the following asymptotic formula holds for the zeros $\theta_n$ of the function $\sigma(\lambda)$:

$$\theta_n = n - \frac{1}{2}\text{sign}\,n - \frac{B_3}{n} + \frac{\eta_n}{n},\quad \{\eta_n\} \in l_2. \qquad (30)$$

Let $\beta = (1+\alpha^2)(B - B_3 \pi)$. Using (28), for the function

$$v_1(\lambda) = \delta(\lambda) - 2\omega - (\alpha\lambda + \beta)\sigma(\lambda) \qquad (31)$$

we obtain the following representation:

$$v_1(\lambda) = \lambda \sin \lambda \pi + (B_1 - \beta)\cos \lambda \pi + \int_{-\pi}^{\pi} p(t)e^{i\lambda t}\,dt,\quad p(t) \in L_2[-\pi, \pi]. \qquad (32)$$

As in [37], define the function

$$v_2(\lambda) = v_1(\lambda) - 2\omega^2 \left[ \frac{\sin \lambda \pi}{\lambda} + B_4 \pi \frac{4 \cos \lambda \pi}{4\lambda^2 - 1} + \frac{r(\lambda)}{\lambda^2} \right], \tag{33}$$

which satisfies the condition

$$v_2(\theta_n) = (-1)^{n+1} \sigma_n \sqrt{v_1^2(\theta_n) - 4\omega^2}, \tag{34}$$

where $B_4$ is some constant, $r(\lambda) = 2\sigma(\lambda) \sum_{n=1}^{\infty} \frac{\theta_n r_n}{\sigma'(\theta_n)(\lambda^2 - \theta_n^2)}$ is an even entire function of exponential type not greater than $\pi$ belonging to $L_2(-\infty, \infty)$ and $\{r_n\} \in l_2$. In view of (32) and (33), for the function

$$\lambda s(\lambda) = \frac{\lambda}{2\omega^2} [v_1(\lambda) - v_2(\lambda)] \tag{35}$$

we get the representation

$$\lambda s(\lambda) = \sin \lambda \pi + B_4 \pi \frac{4\lambda \cos \lambda \pi}{4\lambda^2 - 1} + \frac{r(\lambda)}{\lambda}. \tag{36}$$

Therefore, by Lemma 3.4.2 of [17], if we denote by $\lambda_n$ $(n = \pm 1, \pm 2, ...)$ the zeros of the function $\lambda s(\lambda)$, then the following asymptotic formula will hold true for them:

$$\lambda_n = n - \frac{B_4}{n} + \frac{\xi_n}{n}, \quad \{\xi_n\} \in l_2. \tag{37}$$

As $\sigma(\theta_n) = 0$, it follows from (31) that $v_1(\theta_n) = \delta(\theta_n) - 2\omega$. Assuming $\lambda = 0$ in (27), we get

$$\delta(0) = \pi \sqrt{1 + \alpha^2} \mu_{-0} \cdot \mu_{+0} \prod_{\substack{k=-\infty \\ k \neq 0}}^{\infty} \frac{\mu_k}{k}. \tag{38}$$

As $\operatorname{sign} \mu_k = \operatorname{sign} k$ for $k = \pm 1, \pm 2, ...$ and $\mu_{-0} < 0$, $\mu_{+0} > 0$ (by condition 2), it follows from (38) that $\delta(0) < 0$. Then, according to the second and the third conditions of the theorem, the inequalities $v_1(\theta_{2n-1}) \geq 2|\omega|$, $v_1(\theta_{2n}) \leq -2|\omega|$ are true. Therefore there exists a number $h_n$ such that

$$v_1(\theta_n) = 2|\omega|(-1)^{n+1} \operatorname{ch} h_n. \tag{39}$$

From (34) and (39) we obtain

$$v_2(\theta_n) = 2|\omega|(-1)^{n+1} \sigma_n |\operatorname{sh} h_n|. \tag{40}$$

Letting $\lambda = \theta_n$ in (35) and taking into account (39), (40), we obtain

$$s(\theta_n) = \frac{1}{2\omega^2}[v_1(\theta_n) - v_2(\theta_n)] = \frac{(-1)^{n+1}}{|\omega|}(\operatorname{ch} h_n - \sigma_n |\operatorname{sh} h_n|) =$$

$$= \frac{(-1)^{n+1} \operatorname{ch} h_n}{|\omega|}(1 - \sigma_n |\operatorname{th} h_n|).$$

Hence, from the obvious inequality $|\operatorname{th} h_n| < 1$ it follows that

$$\operatorname{sign} s(\theta_n) = (-1)^{n+1}. \tag{41}$$

Then, according to the second condition of the theorem, each of intervals

$$..., (\theta_{-3}, \theta_{-2}), (\theta_{-2}, \theta_{-1}), (\theta_1, \theta_2), (\theta_2, \theta_3), ...$$

contains one and, by the asymptotic formula (37), only one zero of the function $s(\lambda)$. Consequently, the zeros $..., \theta_{-2}, \theta_{-1}, \theta_1, \theta_2, ...$ of the function $\sigma(\lambda)$ and the zeros $..., \lambda_{-2}, \lambda_{-1}, \lambda_1, \lambda_2, ...$ of the function $s(\lambda)$ satisfy the inequalities

$$... < \theta_{-3} < \lambda_{-2} < \theta_{-2} < \lambda_{-1} < \theta_{-1} < 0 < \theta_1 < \lambda_1 < \theta_2 < \lambda_2 < \theta_3 < .... \tag{42}$$

Let at last

$$s_1(\lambda) = \sigma(\lambda) - \gamma s(\lambda), \tag{43}$$

where $\gamma = (B_4 - B_3)\pi$. The equality (43), together with the relations (29) and (36), implies the representation

$$s_1(\lambda) = \cos \lambda\pi - B_4 \pi \frac{\sin \lambda\pi}{\lambda} + \int_0^\pi r_1(t) \frac{\sin \lambda t}{\lambda} dt$$

(where $r_1(t) \in L_2[0, \pi]$), which, in turn, due to Lemma 3.4.2 of [17], implies the following asymptotic formula for the zeros $v_n$ ($v_{-n} = -v_n$, $n = \pm 1, \pm 2, ...$) of the function $s_1(\lambda)$:

$$v_n = n - \frac{1}{2}\operatorname{sign} n - \frac{B_4}{n} + \frac{\tau_n}{n}, \quad \{\tau_n\} \in l_2. \tag{44}$$

On the other hand, by (41)-(43) we have $\operatorname{sign} s_1(\lambda_n) = \operatorname{sign} \sigma(\lambda_n) = (-1)^n$. Then it is easy to see that $v_m^2 < \lambda_m^2 < v_{m+1}^2$ ($m = 1, 2,...$). Thus, the zeros of the functions $\sqrt{\lambda} s(\sqrt{\lambda})$ and $s_1(\sqrt{\lambda})$ alternate and satisfy the asymptotic formulas (37) and (44). Hence it follows (see Theorem 3.4.1 of [17]) that there exists a unique real function $q(x) \in L_2[0, \pi]$ such that the considered sequences of zeros are the spectra of the boundary value problems generated on the interval $[0, \pi]$ by the same equation (1) with the found coefficient $q(x)$ and the boundary conditions (22), (25), and the equalities $s(\lambda)=s(\pi,\lambda)$, $s_1(\lambda)=s'(\pi,\lambda)$ hold. Taking into account these equalities, it is easy to prove that the spectrum of the constructed boundary value problem coincides with the sequence $\{\mu_k\}$. Theorem is proved.

It can be shown that if the operator $A$ given by the equalities
$$Ay = -y'' + q(x)y,$$
$$D(A) = \{y \in W_2^2[0,\pi] : y'(0) + \beta y(0) + \omega y(\pi) = y'(\pi) + \gamma\, y(\pi) - \omega y(0) = 0\}$$
is positive, then the conditions of Theorem 2 are also necessary.

**Acknowledgments.** This work was supported by the Science Development Foundation under the President of the Republic of Azerbaijan –
Grant № EİF/MQM/Elm-Tehsil-1-2016-1(26)-71/05/1.